# ON THE CONVERGENCE FROM DISCRETE TO CONTINUOUS TIME IN AN OPTIMAL STOPPING PROBLEM[1]

By Paul Dupuis and Hui Wang

*Brown University*

We consider the problem of optimal stopping for a one-dimensional diffusion process. Two classes of admissible stopping times are considered. The first class consists of all nonanticipating stopping times that take values in $[0, \infty]$, while the second class further restricts the set of allowed values to the discrete grid $\{nh : n = 0, 1, 2, \ldots, \infty\}$ for some parameter $h > 0$. The value functions for the two problems are denoted by $V(x)$ and $V^h(x)$, respectively. We identify the rate of convergence of $V^h(x)$ to $V(x)$ and the rate of convergence of the stopping regions, and provide simple formulas for the rate coefficients.

**1. Introduction.** One of the classical formulations of stochastic optimal control is that of optimal stopping, where the only decision to be made is when to stop the process. Upon stopping, a benefit is received (or a cost is paid), and the objective is to maximize the expected benefit (or minimize the expected cost). Although the formulation is very simple, this optimization problem has many practical applications. Examples include the pricing problems in investment theory, the valuation of American options, the development of natural resources and so on [1, 2, 4, 5, 8, 11, 12, 19, 20, 21, 22].

The formulation of the optimal stopping problem requires the specification of the class of allowed stopping times. Typically, one assumes these to be nonanticipative in an appropriate sense, so that the control does not have knowledge of the future. Another important restriction is with regard to the actual time values at which one can stop, and here there are two important cases: continuous time and discrete time. In the first case, the stopping time is allowed to take values in the interval $[0, \infty]$, with $\infty$ corresponding to the

Received June 2001; revised May 2004.

[1]Supported in part by NSF Grants DMS-00-72004, DMS-01-03669 and Army Research Office Grant ARO-DAAD19-99-1-0223.

decision to never stop. In the second case, there is a fixed discrete set of times, and the stopping time must be selected from this set. Typically, this discrete set is a regular grid, for example, $D^h \doteq \{nh : n \in \mathbb{N}_0 \cup \{\infty\}\}$, where $h > 0$ is the grid spacing.

In the present paper we focus exclusively on the one-dimensional case. Although a statement of precise assumptions is deferred to Section 2, a rough description of the continuous and discrete time problems is as follows.

*Continuous time optimal stopping.* We use the stochastic process model

$$\frac{dS_t}{S_t} = b(S_t)\,dt + \sigma(S_t)\,dB_t,$$

where $b$ and $\sigma$ are bounded continuous functions from $\mathbb{R}$ to $\mathbb{R}$, and $B$ is a standard Brownian motion. Although the results can be extended to cover other diffusion models as well, we focus on this model because of its wide use in optimal stopping problems that occur in economics and finance. We consider a payoff defined in terms of a nondecreasing function $\phi : \mathbb{R} \to [0, \infty)$. The payoff from stopping at time $t$ is $\phi(S_t)$, and the decision maker wants to maximize the expected present value by judiciously choosing the stopping time. This is modeled by the optimal stopping problem with value function

$$V(x) \doteq \sup_{\tau \in \mathcal{S}} E[e^{-r\tau}\phi(S_\tau)|S_0 = x],$$

where $r > 0$ is the discount rate and $\mathcal{S}$ is the set of all admissible stopping times. The stopping times are allowed to take values in $[0, \infty]$. Let

$$\mathcal{L}V(x) = \tfrac{1}{2}\sigma^2(x)x^2 V''(x) + b(x)x V'(x).$$

Then the dynamic programming equation for this problem is

$$\max[\phi(x) - V(x), \mathcal{L}V(x) - rV(x)] = 0.$$

If $\phi$ is convex and nondecreasing, it is often optimal to stop when the process $S_t$ first exceeds some fixed threshold $x_*$. In this case, the value function $V(x)$ equals $\phi(x)$ for $x \geq x_*$, and it satisfies the ordinary differential equation $-rV(x) + \mathcal{L}V(x) = 0$ for $x < x_*$.

*Discrete time optimal stopping.* The process model is the same as before, but the set of possible stopping times is restricted to those that take values in the time grid $D^h \doteq \{nh : n \in \mathbb{N}_0 \cup \{\infty\}\}$. The optimal strategy is often similar to the continuous time case: stop the first time $S_{nh}$ exceeds some fixed threshold $x_*^h$. Let $V^h(x)$ denote the value function. The pair $(V^h(x), x_*^h)$ satisfy the dynamic programming equation [25]

$$V^h(x) = \begin{cases} \phi(x), & x \in [x_*^h, \infty), \\ e^{-rh}E[V^h(S_h)|S_0 = x], & x \in (0, x_*^h). \end{cases}$$

Closed-form solutions to this dynamic programming equation are not usually available.

The aim of the present paper is to examine the connection between these two optimal stopping problems as $h \to 0$. There are two questions of main interest:

- What is the convergence rate of the optimal exercise boundary $x_*^h$ to $x_*$, and what is the rate coefficient?
- What is the convergence rate of the value function $V^h(x)$ to $V(x)$, and what is the rate coefficient?

The goal is to use the more readily available solution to the continuous time problem to approximate the solution in discrete time. As we will see in Section 2, the optimal exercise boundaries converge with rate $\sqrt{h}$, while the value functions converge with rate $h$. In both cases there is a well-defined rate coefficient. The coefficient in the case of the exercise boundary is defined in terms of the expected value of a functional of local time of Brownian motion, while the coefficient for the value function involves both local time and excursions of Brownian motion.

Few existing results are concerned with the rate of convergence of approximations for this class of problems. Lamberton [17, 18] considers the binomial tree approximation for pricing American options and its generalizations in order to obtain upper and lower bounds (though not a rate of convergence) for the value function. The pricing of American options is equivalent to solving a finite-horizon optimal stopping problem, and there is no closed form solution. The goal in [17, 18] is, in fact, opposite that of the present paper, in that the discrete time problem is used to approximate the continuous time problem.

As we have noted previously, the motivation for this study is to exploit situations where the continuous time problem can be more or less solved explicitly (e.g., the one-dimensional problems considered in the present work). Our results allow one to explicitly compute accurate approximations for the discrete time problem, and thus avoid numerical approximation. Whether or not one can find a precise rate of convergence and rate coefficients for the analogous question in numerical approximation (where both time and state are discretized) is an interesting open question.

The outline of the paper is as follows. In Section 2 we introduce notation and define the basic optimization problems. Two important universal constants are introduced in Section 3. In Section 4 we state the main result, give an illustrative example, and then lay out the main steps in the proof of the approximation theorem. The proofs of two key approximations which are intimately connected with the local time and excursions of Brownian motion are given in Section 5. The paper concludes with an Appendix in which (i) a result on a conditional distribution of the exit time is proved, and (ii) representations for the universal constants are derived.

**2. Notation, assumptions and background.** Consider a probability space $(\Omega, \mathcal{F}, P; \mathbb{F})$ with filtration $\mathbb{F} = (\mathcal{F}_t)$ satisfying the usual conditions: right-continuity and completion by $P$-negligible sets. The *state process* $S = (S_t, \mathcal{F}_t)$ is modeled by

$$\frac{dS_t}{S_t} = b(S_t)\,dt + \sigma(S_t)\,dB_t, \qquad S_0 \equiv x.$$

Here $B = (B_t, \mathcal{F}_t)$ is a standard Brownian motion.

Define the continuous time value function

$$V(x) \doteq \sup_{\tau \in \mathcal{S}} E[e^{-r\tau}\phi(S_\tau)|S_0 = x],$$

where the supremum is over all stopping times with respect to the filtration $\mathbb{F}$. Define the discrete time value function

$$V^h(x) \doteq \sup_{\tau \in \mathcal{S}^h} E[e^{-r\tau}\phi(S_\tau)|S_0 = x],$$

where $\mathcal{S}^h$ is the set of all stopping times that take values in $D^h$.

The following assumptions will be used throughout the paper.

CONDITION 2.1.

1. The coefficients $b \colon \mathbb{R} \to \mathbb{R}$ and $\sigma \colon \mathbb{R} \to \mathbb{R}$ are bounded and continuous, with $\inf_{x \in \mathbb{R}} \sigma(x) > 0$. Furthermore, $xb(x)$ and $x\sigma(x)$ are Lipschitz continuous.
2. $\phi \colon \mathbb{R} \to [0, \infty)$ is nondecreasing, and both $\phi$ and its derivative $\phi'$ are of polynomial growth. Furthermore,
$$\sup_{t \geq 0} e^{-rt}\phi(S_t) \in \mathbb{L}^1, \qquad \lim_{t \to \infty} e^{-rt}\phi(S_t) = 0 \qquad \text{a.s.}$$
3. The "continuation" region for the continuous-time optimal stopping problem takes the form $\{x \colon V(x) > \phi(x)\} = (0, x_*)$.
4. The "continuation" region for the discrete-time optimal stopping problem takes the form $\{x \colon V^h(x) > \phi(x)\} = (0, x_*^h)$.
5. The payoff function $\phi$ is twice continuously differentiable in a neighborhood of $x_*$.
6. The smooth-fit-principle holds, that is, the value function $V$ is $\mathcal{C}^1$ across the optimal exercise boundary $x_*$.

As noted in the Introduction, $V$ satisfies the dynamic programming equation

$$\max[\phi(x) - V(x), \mathcal{L}V(x) - rV(x)] = 0.$$

Note that usually $V$ is only once continuously differentiable across the optimal exercise boundary $x = x_*$. Since $\phi(x) = V(x)$ if $x \in [x_*, \infty)$ and $\phi(x) < V(x)$ if $x \in (0, x_*)$, it follows that $V''(x_*-) \geq \phi''(x_*)$, where the $-$ denotes limit from the left. Define

$$(2.1) \qquad A \doteq \frac{V''(x_*-) - \phi''(x_*)}{\phi(x_*)} \geq 0.$$

Although one can construct examples where $A = 0$, as the next remark shows, the case $A > 0$ is in a certain sense generic. We will assume this condition below, and note that the rate of convergence of the optimal threshold does not depend on $A$ at all.

REMARK 2.1. The change of variable $t = -\log x$ can be used to transform the ordinary differential equation (ODE) $\mathcal{L}f(x) - rf(x) = 0$ on $(0, \infty)$ into the ODE

$$\tfrac{1}{2}\sigma(e^{-t})W''(t) + [\tfrac{1}{2}\sigma(e^{-t}) - b(e^{-t})]W'(t) - rW(t) = 0$$

on $\mathbb{R}$. Since $\sigma(x) > 0$ for $x > 0$, the classical theory for solutions of ODEs [3] can be used to show that the general solution to $\mathcal{L}V(x) - rV(x) = 0$ can be written in the form $c_1 f_1(x) + c_2 f_2(x)$, where $f_1(x)$ is positive and bounded as $x \downarrow 0$ and $f_2(x)$ is unbounded as $x \downarrow 0$. Under Condition 2.1, the function $f_1$ is twice continuously differentiable on $(0, \infty)$. $V(x)$ is then equal to $c_1 f_1(x)$ for $x \in (0, x_*]$ and equal to $\phi(x)$ for $x \in [x_*, \infty)$, where $c_1$ and $x_*$ are determined by the principle of smooth fit, that is,

$$c_1 f_1(x_*) = \phi(x_*) \quad \text{and} \quad c_1 f_1'(x_*) = \phi'(x_*).$$

REMARK 2.2. If $S$ is a geometric Brownian motion with $b(x) \equiv b$ and $\sigma(x) \equiv \sigma$, and $\phi(x) = (x - k)^+$ for some constant $k$, then Condition 2.1 holds when $r > b$ [8]. For $r < b$, the value function for the optimal stopping problem is $+\infty$, and there is no optimal stopping time. For the boundary case $r = b$, the value function $V(x) \equiv x$ and there is no optimal stopping time.

REMARK 2.3. It is usually not a priori clear if parts 3 and 4 of Condition 2.1 hold for a general state process. Counterexamples can be found in [6, 10]. Interested readers may also find the results in [7] helpful.

Below we give a sufficient condition that is very easy to verify in the case $\phi(x) = (x - k)^+$. Suppose parts 1 and 2 of Condition 2.1 hold, that derivatives of $b$ and $\sigma$ exist and are Hölder continuous for some $\delta > 0$, and, in addition, that

$$r \geq \sup_{x \in (0, \infty)} \{b(x) + xb'(x)\}.$$

We claim that parts 3 and 4 of Condition 2.1 hold. We will show that part 3 holds and omit the analogous proof for 4. Now fix $x \geq y$. We have [15]

$$Z_t \doteq S_t^x - S_t^y = \int_y^x D_t^z \, dz,$$

where $D_t^z = \partial S_t^z / \partial z$ satisfies the SDE

$$\frac{dD_t^z}{D_t^z} = [b(S_t^z) + S_t^z b'(S_t^z)] \, dt + [\sigma(S_t^z) + S_t^z \sigma'(S_t^z)] \, dB_t, \qquad D_0^z = 1.$$

Note that $D^z$ and $Z$ are both nonnegative processes. The condition on $r$ implies $\{e^{-rt} D_t^z\}$ is a supermartingale, and, therefore, so is $\{e^{-rt} Z_t\}$. Observe that $x \geq y$ implies

$$\phi(x) - \phi(y) = (x - k)^+ - (y - k)^+ \leq x - y.$$

Thus, for any stopping time $\tau \in [0, \infty]$,

$$E e^{-r\tau} [\phi(S_\tau^x) - \phi(S_\tau^y)] \leq E[e^{-r\tau} Z_\tau] \leq Z_0 = x - y,$$

where the second inequality follows from the optional sampling theorem. It follows immediately that, for all $x \geq y$,

$$V(x) - V(y) \leq x - y.$$

This implies that

$$\{x : V(x) = \phi(x)\} = [x^*, \infty)$$

for some real number $x^*$. Indeed, if $V(y) = \phi(y) = (y - k)^+$, then since $V > 0$, we must have $y > k$. It follows that, for all $x \geq y$,

$$V(x) \leq V(y) + (x - y) = (y - k) + (x - y) = x - k = \phi(x).$$

But $V \geq \phi$ trivially, whence $V(x) = \phi(x)$ for all $x \geq y$. This completes the proof.

REMARK 2.4.   If $S$ is a geometric Brownian motion with $b(x) \equiv b$, and $\sigma(x) \equiv \sigma$, and $\phi(x) = (\sum_i A_i x^{\alpha_i} - k)^+$ for some positive constants $(A_i, \alpha_i)$ and $k \geq 0$, then one can show that $V(x) - \phi(x)$ is decreasing, which in turn implies that parts 3 and 4 of Condition 2.1 hold. A similar argument can be found in [11, 14].

REMARK 2.5.   We wish to point out that part 6 of Condition 2.1 (i.e., the principle of smooth fit) is not an ad hoc assumption. Much research has been done on the validity of this principle under various conditions, especially for the one-dimensional diffusion case. Interested readers may find the list of references [6, 11, 14, 24] useful.

**3. Two universal constants.** In this section we introduce a pair of universal constants that play an important role in determining the rate coefficients of the convergence.

Let $B$ be a standard one-dimensional Brownian motion. For a fixed constant $u \in [0, 1)$, define the process $W = \{W_t, t \geq u\}$ by

$$W_t \doteq B_t - B_u.$$

In other words, $W$ is a Brownian motion starting at time $t = u$ with initial condition 0. Let $N \doteq \inf\{n \in \mathbb{N} : W_n \geq 0\}$. Note that $N$ is finite with probability one. Define

$$(3.1) \qquad H(u) \doteq EW_N^2 \quad \text{and} \quad M(u) \doteq EW_N.$$

In terms of these functions, we define the constants

$$(3.2) \qquad \Theta = \int_0^1 H(u)\,du \quad \text{and} \quad \Gamma = \int_0^1 M(u)\,du.$$

Note that $H(u) > M^2(u)$, and, therefore,

$$\Theta = \int_0^1 H(u)\,du > \int_0^1 M^2(u)\,du \geq \left(\int_0^1 M(u)\,du\right)^2 = \Gamma^2.$$

Expressions (3.1) and (3.2) for $\Theta$ and $\Gamma$ are useful for purposes of approximation (e.g., Monte Carlo simulation). However, the next lemma connects them with the quantities that will actually arise in the approximation of the optimal stopping problem. The proof of the lemma is given in the Appendix.

LEMMA 3.1. *For every fixed $u \in [0, 1)$,*

$$H(u) = E \int_u^N \mathbb{1}_{\{W_t \geq 0\}}\,dt \quad \text{and} \quad M(u) = E L_{u,N}^W(0),$$

*where $L_{u,N}^W(0)$ is the local time of $W$ on the interval $[u, N]$.*

REMARK 3.1. We have employed Monte Carlo simulation to obtain the approximations $\Theta \approx 0.589$ and $\Gamma \approx 0.582$.

**4. The approximation theorem.** Our main result is the following.

THEOREM 4.1. *Assume Condition 2.1, and define the constants $A$, $\Theta$ and $\Gamma$ by (2.1) and (3.2). Assume that $A > 0$. The following conclusions hold for all $x \in (0, x_*)$:*

1.

$$\frac{V^h(x) - V(x)}{V(x)} = -\frac{1}{2} A x_*^2 \sigma^2(x_*)(\Theta - \Gamma^2)h + o(h).$$

2.

$$x_*^h = x_* - \Gamma x_* \sigma(x_*) \sqrt{h} + o(\sqrt{h}).$$

EXAMPLE 4.1.   Consider the special case where $b(x) \equiv b$ and $\sigma(x) \equiv \sigma$. Assume $r > b$ and $\phi(x) = (x - k)^+$ for some constant $k > 0$. It follows that the value function for the continuous time optimal stopping problem is

$$V(x) = \begin{cases} Bx^\alpha, & x < x_*, \\ x - k, & x \geq x_*, \end{cases}$$

where

$$\alpha = \left( \frac{1}{2} - \frac{b}{\sigma^2} \right) + \sqrt{ \left( \frac{1}{2} - \frac{b}{\sigma^2} \right)^2 + \frac{2r}{\sigma^2} }, \qquad B = \frac{x_* - k}{x_*^\alpha}$$

and

$$x_* = \frac{\alpha}{\alpha - 1} k.$$

According to the theorem,

$$x_*^h = x_*(1 - \Gamma \sigma \sqrt{h}) + o(\sqrt{h}) = \frac{\alpha k}{\alpha - 1}(1 - \Gamma \sigma \sqrt{h}) + o(\sqrt{h})$$

and

$$A = \frac{V''(x_*-) - \phi''(x_*)}{x_* - k} = \frac{B\alpha(\alpha - 1)(x_*)^{\alpha - 2}}{x_* - k},$$

which, after some algebra, yields

$$\frac{V^h(x) - V(x)}{V(x)} = -\frac{1}{2}\alpha(\alpha - 1)(\Theta - \Gamma^2)\sigma^2 h + o(h).$$

4.1.  *Overview of the proof.*   In this section we outline and prove the main steps in the proof of Theorem 4.1. The proofs of two key asymptotic expansions are deferred to the next section.

To simplify the analysis, we first introduce a bounded modification of the payoff function $\phi$. This modification will not affect the asymptotics at all; see Proposition 4.2.

Let $\bar{\phi} \leq \phi$ be an increasing function satisfying

(4.1) $$\bar{\phi}(x) = \begin{cases} \phi(x), & \text{if } x \leq x_* + a, \\ k, & \text{if } x \geq x_* + 2a. \end{cases}$$

Here $a$ and $k$ are two positive constants, whose specific values are not important. Without loss of generality, we assume that $\bar{\phi}$ is twice continuously differentiable in the region $[x_*, \infty)$. Suppose $h$ and $\delta$ are two positive constants, and let $x_\delta \doteq x_* - \delta$. We consider the quantities

$$\bar{W}_\delta(x) \doteq E^x[e^{-r\tau_\delta}\bar{\phi}(S_{\tau_\delta})] \quad \text{and} \quad W_\delta(x) = E^x[e^{-r\tau_\delta}\phi(S_{\tau_\delta})],$$

where

$$\tau_\delta \doteq \inf\{t \geq 0 : S_t \geq x_\delta\},$$

and $E^x$ denotes expectation conditioned on $S_0 = x$. Note that $W_\delta(x) = \bar{W}_\delta(x)$ for all $x \leq x_\delta$. We also define

$$\bar{W}_\delta^h(x) \doteq E^x[e^{-r\tau_\delta^h}\bar{\phi}(S_{\tau_\delta^h})] \quad \text{and} \quad W_\delta^h(x) \doteq E^x[e^{-r\tau_\delta^h}\phi(S_{\tau_\delta^h})],$$

where

$$\tau_\delta^h \doteq \inf\{nh \geq 0 : S_{nh} \geq x_\delta\}.$$

These are all cost functions for an a priori fixed (and possibly suboptimal) stopping region.

*Main idea of the proof.* The main idea for proving the rates of convergence is as follows. Write

$$W_\delta^h(x) - V(x) = [W_\delta^h(x) - \bar{W}_\delta^h(x)] + [\bar{W}_\delta^h(x) - \bar{W}_\delta(x)] + [\bar{W}_\delta(x) - V(x)].$$

For each term, we will obtain approximations as $h$ and $\delta$ tend to zero. It turns out the leading term for the sum has the following form:

$$-\tfrac{1}{2}a_1\delta^2 + a_2\delta\sqrt{h} + a_3 h + \text{higher-order term}.$$

Here $a_1, a_2$ and $a_3$ are constants with respect to $\delta$ and $h$ (though some will depend on $x$) and with $a_1 > 0$. Since $x_*^h$ is the optimal exercise boundary for the discrete problem, the mapping $\delta \to W_\delta^h(x)$ attains its maximum at $\delta_* = x_* - x_*^h$. Hence, one would expect that $\delta_*$ would approximately maximize the leading term, or

(4.2) $$\delta_* = \frac{a_2}{a_1}\sqrt{h} + o(\sqrt{h}).$$

Substituting this back in, one would further expect

$$V^h(x) = W_{\delta_*}^h(x) = \left(-\frac{a_2^2}{a_1} + a_3\right)h + o(h).$$

Thus, we obtain the precise asymptotic behavior of both the stopping regions and the value functions once the quantities $a_1, a_2$ and $a_3$ are determined.

This is, in fact, how the argument will proceed. We begin with the estimation of the first term, which turns out to be negligible for small $h$ and $\delta$. Define the quantity

(4.3) $$\triangle_{\delta,h} \doteq W_\delta^h(x) - \bar{W}_\delta^h(x) = E^x[e^{-r\tau_\delta^h}(\phi(S_{\tau_\delta^h}) - \bar{\phi}(S_{\tau_\delta^h}))].$$

We have the following result.

PROPOSITION 4.2. *Define* $\triangle_{\delta,h}$ *by* (4.3). *There exist constants* $L < \infty$ *and* $\varepsilon > 0$ *such that*

$$|\triangle_{\delta,h}| \leq Le^{-\varepsilon/h}$$

*for all sufficiently small* $\delta$ *and* $h$.

PROOF. The proof is based on the following bound. Let $a$ be as in the characterization (4.1) of $\bar{\phi}$. Then for any $x \leq x_*$, $y \geq x_* + a$, and $h > 0$, we have

$$(4.4) \qquad P(S_h > y | S_0 = x) \leq \exp\left\{-\frac{1}{c_2 h}\left[\log\frac{y}{x_*} - c_1 h\right]^2\right\},$$

where the positive constants $c_1, c_2$ depend only on the coefficients $b, \sigma$. The proof of this inequality is a standard application of exponential martingales [23], and thus omitted.

We now complete the proof of the proposition. To ease the exposition, we use $\tau$ in lieu of $\tau_\delta^h$ throughout the proof. We have

$$\triangle_{\delta,h} = \sum_{n=1}^\infty e^{-rnh} E^x[(\phi(S_{nh}) - \bar{\phi}(S_{nh}))\mathbb{1}_{\{\tau=n\}}]$$

$$\leq \sum_{n=1}^\infty e^{-rnh} \int_{a+x_*}^\infty |\phi'(y) - \bar{\phi}'(y)| P^x(S_{nh} > y, \tau = n)\, dy.$$

Fix an $n \in \mathbb{N}$. Define the stopping time $\sigma \doteq \inf\{t \geq (n-1)h : S_t \geq x_\delta\}$. Then

$$P^x(S_{nh} > y | \tau = n)$$

$$= \int_0^h P^x(S_{nh} > y | \sigma = nh - t, S_{nh} \geq x_\delta, S_{(n-1)h} < x_\delta, \ldots, S_0 < x_\delta)$$

$$\times P^x(\sigma \in nh - dt | \sigma \leq nh, S_{nh} \geq x_\delta, S_{(n-1)h} < x_\delta, \ldots, S_0 < x_\delta).$$

However, the strong Markov property implies, for all $t \in [0, h]$, that

$$P^x(S_{nh} > y | \sigma = nh - t, S_{nh} \geq x_\delta, S_{(n-1)h} < x_\delta, \ldots, S_0 < x_\delta)$$

$$= P(S_t > y | S_0 = x_\delta, S_t \geq x_\delta)$$

$$= P(S_t > y | S_0 = x_\delta)/P(S_t \geq x_\delta | S_0 = x_\delta).$$

The denominator in this display is uniformly bounded from below away from zero, for all $t \in [0, 1]$:

$$P(S_t \geq x_\delta | S_0 = x_\delta) \geq \alpha > 0.$$

(See the proof of Lemma 5.5 for the detailed calculations in an analogous setting.) Using (4.4), for all small $h > 0$ and $t \in (0, h)$,

$$P(S_t > y | S_0 = x_\delta) \leq \exp\left\{-\frac{1}{c_2 t}\left[\log\frac{y}{x_*} - c_1 t\right]^2\right\}$$

$$\leq \exp\left\{-\frac{1}{c_2 h}\left[\log\frac{y}{x_*} - c_1 h\right]^2\right\}.$$

Now since $\phi'$ is of polynomial growth and $\bar{\phi}'(x)$ is zero for large $x$, it follows that there are finite constants $R$ and $m$ such that

$$|\phi'(y) - \bar{\phi}'(y)| \leq Ry^{m-1} \qquad \text{for all } y > x_* + a.$$

Hence, for all small $\delta > 0$, the change of variable $x = \log(y/x_*) - c_1 h$ gives

$$\triangle_{\delta, h} \leq \frac{R}{\alpha} \sum_{n=1}^{\infty} e^{-rnh} P^x(\tau = n) \int_{a+x_*}^{\infty} y^{m-1} \exp\left\{-\frac{1}{c_2 h}\left[\log\frac{y}{x_*} - c_1 h\right]^2\right\} dy$$

$$= \frac{R}{\alpha}(x_*)^m e^{mc_1 h} \sum_{n=1}^{\infty} e^{-rnh} P^x(\tau = n) \int_{\log(1+a/x_*)-c_1 h}^{\infty} e^{mx - x^2/(c_2 h)} dx.$$

Let $\Phi$ be the cumulative distribution function for the standard normal distribution. For $h$ small enough, there exists positive numbers $\bar{a}, \bar{C}, \bar{c}$ such that

$$\triangle_{\delta, h} \leq \bar{C} \sum_{n=1}^{\infty} e^{-rnh} P^x(\tau = n) \int_{\log(1+\bar{a}/x_*)}^{\infty} e^{-x^2/(\bar{c}h)} dx$$

$$= \bar{C}\sqrt{2\pi\bar{c}} \cdot \Phi\left[-\frac{1}{\sqrt{\bar{c}h}}\log\left(1 + \frac{\bar{a}}{x_*}\right)\right] \cdot \sqrt{h} \sum_{n=1}^{\infty} e^{-rnh} P^x(\tau = n)$$

$$\leq \bar{C}\sqrt{2\pi\bar{c}} \cdot \Phi\left[-\frac{1}{\sqrt{\bar{c}h}}\log\left(1 + \frac{\bar{a}}{x_*}\right)\right] \cdot \sqrt{h}.$$

We complete the proof of the proposition by using the asymptotic relation

$$\Phi(-x) \sim \frac{1}{\sqrt{2\pi}x} e^{-x^2/2}$$

as $x \to \infty$. $\quad\square$

The bound just proved shows that $[W_\delta^h(x) - \bar{W}_\delta^h(x)]$ is exponentially small as $h \to 0$, uniformly for all small $\delta > 0$. We now consider the terms $[\bar{W}_\delta^h(x) - \bar{W}_\delta(x)]$ and $[\bar{W}_\delta(x) - V(x)]$. When considering the asymptotic behavior of these terms, it is often convenient to scale $\delta$ with $h$ as $h \to 0$ in the manner suggested by (4.2). For the remainder of this proof, unless explicitly stated otherwise, we will assume that

$$(4.5) \qquad \delta = c\sqrt{h} + o(\sqrt{h}) \qquad \text{as } h \to 0$$

for a nonnegative parameter $c$. With an abuse of notation, the quantities $\bar{W}^h_\delta(x)$ and $\bar{W}_\delta(x)$ will be denoted by $\bar{W}^h_c(x)$ and $\bar{W}_c(x)$ when the relation (4.5) holds.

We next estimate $[\bar{W}_c(x) - V(x)]$ as $h \to 0$.

PROPOSITION 4.3. *Assume Condition* 2.1 *and define $A$ by* (2.1). *Assume also that $A > 0$. Then*

$$\bar{W}_c(x) - V(x) = [-\tfrac{1}{2}Ac^2h + o(h)]V(x).$$

PROOF. Recall that $V(x)$ can be characterized, for $x \leq x_*$, as a multiple of the bounded (in a neighborhood of zero) solution $f_1$ to $\mathcal{L}f(x) - rf(x) = 0$; see Remark 2.1. $\bar{W}_c(x)$ can be likewise characterized, with the constant determined by the boundary condition $\bar{W}_c(x_\delta) = \bar{\phi}(x_\delta)$. Thus,

$$\bar{W}_c(x) = \frac{\bar{\phi}(x_\delta)}{V(x_\delta)}V(x) \qquad \text{for all } x \in (0, x_\delta].$$

We now apply Taylor's theorem for small $\delta \geq 0$, and use $x_\delta \doteq x_* - \delta, V(x_*) = \bar{\phi}(x_*), V'(x_*) = \bar{\phi}'(x_*)$, and the definition of $A$ to obtain

$$\frac{\bar{W}_c(x) - V(x)}{V(x)} = -\frac{1}{2}A\delta^2 + o(\delta^2).$$

The proof is completed by using (4.5). □

In the next proposition we state the expansion for $[\bar{W}^h_\delta(x) - \bar{W}_\delta(x)]$. This estimate deals with the critical comparison between the discrete and continuous time problems. The proof of this expansion is detailed, and therefore deferred to the next section.

PROPOSITION 4.4. *Assume Condition* 2.1 *and define $A$, $\Theta$ and $\Gamma$ by* (2.1) *and* (3.2). *Assume also that $A > 0$. Then*

$$\bar{W}^h_c(x) - \bar{W}_c(x) = [\Gamma x_*\sigma(x_*)Ach - \tfrac{1}{2}A\Theta x_*^2\sigma^2(x_*)h + o(h)]V(x).$$

PROOF OF THEOREM 4.1. Recall that $x^h_*$ is the optimal boundary for the stopping problem with value function $V^h$. On the stopping region, we always have $V^h(x) = \phi(x)$. Also, since $V^h(x)$ is defined by supremizing over a subset of the stopping times allowed in the definition of $V(x)$, it follows that $V^h(x) \leq V(x)$. Since $V(x) \geq \phi(x)$ for all $x$, it follows that $x^h_* \leq x_*$.

According to Propositions 4.2, 4.3 and 4.4, for each fixed $c \in [0, \infty)$,

$$\frac{W^h_c(x) - V(x)}{V(x)} = \left[ -\frac{1}{2}Ac^2h + \Gamma x_*\sigma(x_*)Ach - \frac{1}{2}A\Theta x_*^2\sigma^2(x_*)h + o(h) \right].$$

This suggests that the choice $c_* \doteq \Gamma x_* \sigma(x_*)$ should define the maximizer and also (at least approximately) the boundary of the optimal stopping region. Inserting this into the last display gives

$$\frac{W_{c_*}^h(x) - V(x)}{V(x)} = \left[\frac{1}{2} A(\Gamma^2 - \Theta) x_*^2 \sigma(x_*)^2 h + o(h)\right],$$

and since $V^h(x) \geq W_{c_*}^h(x)$, it follows that

$$(4.6) \qquad \liminf_{h \downarrow 0} \frac{V^h(x) - V(x)}{V(x)h} \geq \frac{1}{2} A(\Gamma^2 - \Theta) x_*^2 \sigma(x_*)^2.$$

Now define $c^h$ by $x_*^h = x_* - c^h \sqrt{h}$. Since $x_*^h \leq x_*$, we know that $c^h \in [0, \infty)$. By taking a convergent subsequence, we can assume that $c^h \to \bar{c} \in [0, \infty]$. Using an elementary weak convergence argument, one can show that $x_*^h \to x_*$. First assume that $\bar{c} \in (0, \infty)$. If $\bar{c} \neq \Gamma x_* \sigma(x_*)$, then by Propositions 4.2, 4.3 and 4.4, we have

$$\limsup_{h \downarrow 0} \frac{V^h(x) - V(x)}{V(x)h} < \frac{1}{2} A(\Gamma^2 - \Theta) x_*^2 \sigma(x_*)^2,$$

which contradicts (4.6). If $\bar{c} = \infty$, then Propositions 4.2 and 4.4 and an argument analogous to the one used in Proposition 4.3 shows that

$$\frac{V^h(x) - V(x)}{V(x)} = -A(c^h)^2 h[1 + o(1)].$$

Since $(c^h)^2 \to +\infty$, this again contradicts (4.6), and, thus, $\bar{c} = \Gamma x_* \sigma(x_*)$. We extend to the original sequence by the standard argument by contradiction, and Theorem 4.1 follows. $\square$

## 5. Approximations and expansions in terms of local time and excursions of Brownian motion.
In this section we prove Proposition 4.4, which is the expansion $\bar{W}_\delta^h(x) - \bar{W}_\delta(x)$ for small $h > 0$. Throughout this section we assume (4.5), which we repeat here for convenience: $\delta = c\sqrt{h} + o(\sqrt{h})$ as $h \to 0$. Recall also that subscripts of $c$ and $\delta$ may be used interchangeably under this condition. Define the error term

$$\varepsilon(x) = E^x[e^{-rh}\bar{W}_\delta(S_h) - \bar{W}_\delta(x)] \qquad \forall x < x_\delta.$$

Observe that $\bar{W}_\delta$ has the representation

$$\bar{W}_\delta(x) = \begin{cases} \bar{\phi}(x), & \text{for } x \geq x_\delta, \\ -\varepsilon(x) + e^{-rh} E[\bar{W}_\delta(S_h)|S_0 = x], & \text{for } x < x_\delta. \end{cases}$$

It follows from the generalized Itô formula that

$$e^{-rh}\bar{W}_\delta(S_h) - \bar{W}_\delta(x)$$

$$(5.1) \qquad = \int_0^h e^{-rt}[-r\bar\phi(S_t) + \mathcal{L}\bar\phi(S_t)]\mathbb{1}_{\{S_t \geq x_\delta\}}\,dt$$

$$+ \triangle\bar{W}'_\delta(x_\delta)\int_0^h e^{-rt}\,dL_t^S(x_\delta) + \int_0^h e^{-rt}\bar{W}'_\delta(S_t)S_t\sigma(S_t)\,dW_t.$$

Here $L^S$ is the local time for the process $S$, and

$$\triangle\bar{W}'_\delta(x_\delta) \doteq \bar{W}'_\delta(x_\delta+) - \bar{W}'_\delta(x_\delta-).$$

It is straightforward to prove that the stochastic integral has expectation zero, and we arrive at the following result.

LEMMA 5.1. *For every $x \in (0, x_\delta)$,*

$$\varepsilon(x) = E^x\int_0^h e^{-rt}[-r\bar\phi(S_t) + \mathcal{L}\bar\phi(S_t)]\mathbb{1}_{\{S_t \geq x_\delta\}}\,dt$$

$$+ E^x\triangle\bar{W}'_\delta(x_\delta)\int_0^h e^{-rt}\,dL_t^S(x_\delta).$$

Let $\tau_\delta^h \doteq \inf\{nh : S_{nh} \geq x_\delta\}$. Then the discounting and the lemma just stated imply the formula

$$\bar{W}_\delta(x) = -E^x\sum_{n=0}^{\tau_\delta^h/h-1} e^{-rnh}\varepsilon(S_{nh}) + E^x e^{-r\tau_\delta^h}\bar\phi(S_{\tau_\delta^h})$$

$$= -E^x\int_0^{\tau_\delta^h} e^{-rt}[-r\bar\phi(S_t) + \mathcal{L}\bar\phi(S_t)]\mathbb{1}_{\{S_t \geq x_\delta\}}\,dt$$

$$- \triangle\bar{W}'_\delta(x_\delta)E^x\int_0^{\tau_\delta^h} e^{-rt}\,dL_t^S(x_\delta) + E^x e^{-r\tau_\delta^h}\bar\phi(S_{\tau_\delta^h})$$

for all $x < x_\delta$. Observe that the last term is precisely $\bar{W}_\delta^h(x)$ by definition. It follows that

$$\bar{W}_\delta^h(x) - \bar{W}_\delta(x) = E^x\int_0^{\tau_\delta^h} e^{-rt}[-r\bar\phi(S_t) + \mathcal{L}\bar\phi(S_t)]\mathbb{1}_{\{S_t \geq x_\delta\}}\,dt$$

$$+ \triangle\bar{W}'_\delta(x_\delta)E^x\int_0^{\tau_\delta^h} e^{-rt}\,dL_t^S(x_\delta).$$

The proof of Proposition 4.4 is thereby reduced to proving the following two results.

PROPOSITION 5.2. *Assume Condition* 2.1 *and define $A$ and $\Theta$ by* (2.1) *and* (3.2). *Assume also that $A > 0$. Then*

$$E^x \int_0^{\tau_\delta^h} e^{-rt}[-r\bar{\phi}(S_t) + \mathcal{L}\bar{\phi}(S_t)]\mathbb{1}_{\{S_t \geq x_\delta\}}\,dt = [-\tfrac{1}{2}A\Theta x_*^2\sigma^2(x_*)h + o(h)]V(x).$$

PROPOSITION 5.3. *Assume Condition* 2.1 *and define $A$ and $\Gamma$ by* (2.1) *and* (3.2). *Assume also that $A > 0$. Then*

$$(5.2) \qquad \triangle\bar{W}_\delta'(x_\delta)E^x \int_0^{\tau_\delta^h} e^{-rt}\,dL_t^S(x_\delta) = [\Gamma x_*\sigma(x_*)Ach + o(h)]V(x).$$

The proofs of Propositions 5.2 and 5.3 use estimates on the excursions and local time of Brownian motion, respectively, and are given in Sections 5.1 and 5.2.

5.1. *Proof of Proposition* 5.2. We recall Lemma 3.1, which states that

$$H(u) \doteq EW_N^2 = E\int_u^N \mathbb{1}_{\{W_t \geq 0\}}\,dt,$$

for arbitrary $u \in [0, 1)$. Here $W = (B_t - B_u, t \geq u)$ is a Brownian motion starting at time $t = u$ with initial condition $W_u = 0$, and $N \doteq \inf\{n \in \mathbb{N} : W_n \geq 0\}$.

LEMMA 5.4. $H(u)$ *is continuous and bounded on the interval* $[0, 1)$.

PROOF. Define

$$Z_u \doteq \int_u^N \mathbb{1}_{\{W_t \geq 0\}}\,dt.$$

We first show that the family $\{Z_u, u \in [0, 1)\}$ is uniformly integrable [and, in particular, that $H(u)$ is bounded]. Indeed, define

$$c_0 \doteq \int_u^1 \mathbb{1}_{\{W_t \geq 0\}}\,dt, \qquad c_j \doteq \int_j^{j+1} \mathbb{1}_{\{W_t \geq 0\}}\,dt, \qquad j \in \mathbb{N}.$$

The key observation is that if $c_j > 0$, then $W$ must spend some time during the interval $[j, j+1]$ to the right of zero, therefore, the probability that $W_{j+1} > 0$ is at least half. Thus, for all $j \in \mathbb{N}_0$,

$$P(N = j+1 \mid N > j, c_j > 0) \geq \tfrac{1}{2}.$$

Let $X_u \doteq \sum_{j=0}^{N-1} \mathbb{1}_{\{c_j > 0\}}$. Clearly, $X_u$ dominates $Z_u$. Furthermore, the strong Markov property implies that

$$P(X_u \geq j+1 \mid X_u \geq j) \leq 1 - \tfrac{1}{2} = \tfrac{1}{2}.$$

This, in turn, implies that $P(X_u \geq n + 1) \leq 2^{-n}$, and, thus,

$$E(Z_u^2) \leq E(X_u^2) = \sum_{n=1}^{\infty} 2nP(X_u \geq n) \leq \sum_{n=1}^{\infty} \frac{n}{2^{n-2}} < \infty.$$

Therefore, $\{Z_u, u \in [0, 1)\}$ is uniformly integrable.

As for the continuity, we write

$$H(u) = E \int_u^N \mathbb{1}_{\{B_t - B_u \geq 0\}} \, dt = EZ_u.$$

Fix any $u \in [0, 1)$ and let $\{u_n\}$ be an arbitrary sequence in $[0, 1)$ with $u_n \to u$. Since for $P(B_n - B_u = 0) = 0$ for every fixed $n$, we have $Z_{u_n} \to Z_u$ with probability one. Since the $Z_{u_n}$ are uniformly integrable, $H(u_n) \to H(u)$. This completes the proof. $\square$

Now for any $u \in [0, 1)$ and $h > 0$, define the function

$$F(h; u) \doteq \frac{1}{h} E \int_{uh}^{N^h h} e^{-r(t-uh)}[-r\bar{\phi}(S_t) + \mathcal{L}\bar{\phi}(S_t)]\mathbb{1}_{\{S_t \geq x_\delta\}} \, dt,$$

where

$$\frac{dS_t}{S_t} = b(S_t) \, dt + \sigma(S_t) \, dB_t, \qquad S_{uh} = x_\delta$$

and

$$N^h \doteq \inf\{n \in \mathbb{N} : S_{nh} \geq x_\delta\}.$$

Let $\lfloor a \rfloor$ denote the integer part of $a$. It follows from strong Markov property that

(5.3)
$$E^x \int_0^{\tau_\delta^h} e^{-rt}[-r\bar{\phi}(S_t) + \mathcal{L}\bar{\phi}(S_t)]\mathbb{1}_{\{S_t \geq x_\delta\}} \, dt$$
$$= hE^x\left[e^{-r\tau_\delta} F\left(h; \frac{\tau_\delta}{h} - \left\lfloor \frac{\tau_\delta}{h} \right\rfloor\right)\right].$$

Consider the change of variable $t \mapsto th$ and the transformation

$$Y_t^{(h)} \doteq \frac{S_{th} - x_\delta}{\sqrt{h}}.$$

We can rewrite

$$F(h; u) = E \int_u^{N^h} e^{-r(t-u)h}[-r\bar{\phi} + \mathcal{L}\bar{\phi}](\sqrt{h}Y_t^{(h)} + x_\delta)\mathbb{1}_{\{Y_t^{(h)} \geq 0\}} \, dt,$$

where $Y^{(h)}$ follows the dynamics

$$dY_t^{(h)} = (\sqrt{h}Y_t^{(h)} + x_\delta)[\sqrt{h}b(\sqrt{h}Y_t^{(h)} + x_\delta) \, dt + \sigma(\sqrt{h}Y_t^{(h)} + x_\delta) \, dB_t]$$

with initial condition $Y_u^{(h)} = 0$.

We have the following result regarding $F(h; u)$. Although part of the proof is similar to that of Lemma 5.4, we provide the details for completeness.

LEMMA 5.5.   1. $F(h; u)$ *is uniformly bounded for small $h$ and $u \in [0, 1)$.*
2.

$$\lim_{h \to 0} F(h; u) = [-r\bar\phi + \mathcal{L}\bar\phi](x_*)H(u),$$

*and the convergence is uniform on any compact subset of* $[0, 1)$.

PROOF.   Consider the family of random variables $\{Z_{h,u} : u \in [0, 1), h \in (0, 1)\}$, where

$$Z_{h,u} \doteq \int_u^{N^h} e^{-r(t-u)h}[-r\bar\phi + \mathcal{L}\bar\phi](\sqrt{h}Y_t^{(h)} + x_\delta)\mathbb{1}_{\{Y_t^{(h)} \geq 0\}} \, dt.$$

We first show this family is uniformly integrable. Since $(-r\bar\phi + \mathcal{L}\bar\phi)$ is bounded, it is sufficient to show that

$$(5.4) \qquad\qquad X_{h,u} \doteq \int_u^{N^h} \mathbb{1}_{\{Y_t^{(h)} \geq 0\}} \, dt$$

are uniformly integrable. Define

$$c_0^{(h)} \doteq \int_u^1 \mathbb{1}_{\{Y_t^{(h)} \geq 0\}} \, dt, \qquad c_j^{(h)} \doteq \int_j^{j+1} \mathbb{1}_{\{Y_t^{(h)} \geq 0\}} \, dt, \qquad j \in \mathbb{N}.$$

As in the proof of Lemma 5.4, if $c_j^{(h)} > 0$, then $Y_t^{(h)}$ spends some time to the right of zero in the interval $[j, j + 1]$. We claim that the probability of $Y_{j+1}^{(h)} \geq 0$ is bounded from below by a positive constant:

$$P(N^h = j + 1 | N^h > j, c_j^{(h)} > 0) \geq \alpha > 0 \qquad \forall \, u \in [0, 1), h \in (0, 1).$$

To this end, it suffices to show that, for some $\alpha > 0$,

$$p_{t,h} \doteq P(Y_t^{(h)} \geq 0 | Y_0^{(h)} \geq 0) \geq \alpha > 0 \qquad \forall \, t \in [0, 1].$$

However, it is easy to see that

$$\begin{aligned}
p_{t,h} &= P(S_{th} \geq x_\delta | S_0 \geq x_\delta) \\
&\geq P\Big(\exp\Big\{\int_0^{th} [b(S_u) - \tfrac{1}{2}\sigma^2(S_u)] \, du + \int_0^{th} \sigma(S_u) \, dB_u\Big\} \geq 1\Big) \\
&\geq P\Big(\int_0^{th} \sigma(S_u) \, dB_u \geq c_1 th\Big),
\end{aligned}$$

where $c_1 \doteq \|b\|_\infty + \|\sigma^2\|_\infty/2$. We can view the stochastic integral $Q_t \doteq \int_0^t \sigma(S_u)\, dB_u$ a time-changed Brownian motion. Indeed, there exists a Brownian motion $W$ such that [13]

$$Q_t = W_{\langle Q \rangle_t}.$$

Let

$$\underline{\sigma} \doteq \inf_x \sigma(x), \qquad \bar{\sigma} \doteq \sup_x \sigma(x).$$

Then

$$\underline{\sigma}^2 h \leq \langle Q \rangle_t \leq \bar{\sigma}^2 h.$$

It follows that

$$p_{t,h} \geq P\left(\min_{\underline{\sigma}^2 th \leq s \leq \bar{\sigma}^2 th} W_s \geq c_1 th\right) = P\left(\min_{\underline{\sigma}^2 \leq s \leq \bar{\sigma}^2} W_s \geq c_1 \sqrt{th}\right),$$

where the last equality follows since $\{W_{ths}/\sqrt{th}, s \geq 0\}$ is still a standard Brownian motion. For $h \in (0,1)$ and $t \in [0,1]$, we can choose

$$\alpha = P\left(\min_{\underline{\sigma}^2 \leq s \leq \bar{\sigma}^2} W_s \geq c_1\right) > 0,$$

which will serve as a lower bound.

Now define

$$M_{h,u} \doteq \sum_{j=0}^{N^h - 1} \mathbb{1}_{\{c_j^{(h)} > 0\}},$$

which clearly dominates $X_{h,u}$. By the strong Markov property,

$$P(M_{h,u} > j+1 | M_{h,u} > j) \leq 1 - \alpha,$$

and, thus,

$$P(M_{h,u} \geq j) \leq (1-\alpha)^{j-1}.$$

This implies that

$$E(M_{h,u}^2) = \sum_{j=0}^\infty 2j P(M_{h,u} \geq j) \leq \sum_{j=0}^\infty 2j(1-\alpha)^{j-1} < \infty,$$

which implies the uniform integrability of $\{Z_{h,u}, u \in [0,1), h \in (0,1)\}$. In particular, $F(h; u)$ is uniformly bounded for $u \in [0,1)$ and $h \in (0,1)$.

For the uniform convergence, it suffices to show that, for any $u \in [0,1)$ and any sequence $u^h \in [0,1)$ converging to $u$,

$$F(h; u^h) = E Z_{h, u^h} \to [-r\bar{\phi} + \mathcal{L}\bar{\phi}](x_*) H(u).$$

Let $Y^{(h)}$ be the process with $Y_{u^h}^{(h)} = 0$. As $h \to 0$, we have that $Y^{(h)}$ converges weakly to $Y$, where $Y$ is defined as

$$Y_t = x_* \sigma(x_*)(B_t - B_u).$$

By the Skorohod representation, we can assume $Y^{(h)} \to Y$ with probability one. Using the uniform integrability, it suffices to show that

$$Z_{h,u^h} \to Z \doteq [-r\bar{\phi} + \mathcal{L}\bar{\phi}](x_*) \int_u^N \mathbb{1}_{\{Y_t \geq 0\}} \, dt$$

with probability one. Note that $N$ is almost surely finite, and that $N^h \to N$ with probability one. The almost sure convergence of $Z_{h,u^h}$ to $Z$ then follows from the dominated convergence theorem, which completes the proof. $\square$

Returning to the proof of Proposition 5.2, we claim that

$$(5.5) \quad \lim_{h \to 0} E^x \left[ e^{-r\tau_\delta} F\left( h; \frac{\tau_\delta}{h} - \left\lfloor \frac{\tau_\delta}{h} \right\rfloor \right) \right] = \Theta[-r\bar{\phi} + \mathcal{L}\bar{\phi}](x_*) E^x[e^{-r\tau_*}].$$

To ease notation, let

$$U_h \doteq \frac{\tau_\delta}{h} - \left\lfloor \frac{\tau_\delta}{h} \right\rfloor.$$

It suffices to show that

$$\lim_{h \to 0} E^x[e^{-r\tau_\delta} F(h; U_h)] = \Theta[-r\bar{\phi} + \mathcal{L}\bar{\phi}](x_*) E^x[e^{-r\tau_*}].$$

In Proposition A.3 in the Appendix we show the following (not very surprising) result. As $h$ and $\delta$ tend to zero, $(\tau_\delta, U_h)$ converges in distribution to $(\tau_*, U)$, where $U$ is uniformly distributed and independent of $\tau_*$. More precisely, we have

$$E^x[e^{-r\tau_\delta} H(U_h)] \to E^x[e^{-r\tau_*}] \int_0^1 H(u) \, du = \Theta E^x[e^{-r\tau_*}].$$

Therefore, to prove (5.5), we must show that

$$\triangle \doteq E^x[e^{-r\tau_\delta} F(h; U_h)] - [-r\bar{\phi} + \mathcal{L}\bar{\phi}](x_*) E^x[e^{-r\tau_\delta} H(U_h)] \to 0.$$

Due to the uniform boundedness of $F$ and $H$, there exists $R \in (0, \infty)$ such that

$$|F(h, u)| + |[-r\bar{\phi} + \mathcal{L}\bar{\phi}](x_*) H(u)| \leq R \qquad \forall u \in [0, 1),$$

when $h$ is small enough. Since $U_h \Rightarrow U$, for $h$ small enough,

$$P(U_h > 1 - \varepsilon) \leq 2\varepsilon.$$

Also, by Lemma 5.5 for $h$ small enough,

$$\sup_{u \in [0, 1-\varepsilon]} |F(h, u) - [-r\bar{\phi} + \mathcal{L}\bar{\phi}](x_*)H(u)| \le \varepsilon.$$

It follows that, for $h$ small enough,

$$\triangle \le \varepsilon P(U_h \le 1 - \varepsilon) + RP(U_h > 1 - \varepsilon) \le (2R+1)\varepsilon,$$

which completes the proof of (5.5).

It follows directly from the definitions of $V(x)$ and $\tau_*$ that

(5.6)                           $E^x[e^{-r\tau_*}] = V(x)/\bar{\phi}(x_*).$

Also, the definition of $A$ in (2.1) and the fact that $(-rV + \mathcal{L}V)(x_*-) = 0$ imply that

$$\begin{aligned}
(-r\bar{\phi} &+ \mathcal{L}\bar{\phi})(x_*) \\
&= (-rV + \mathcal{L}V)(x_*-) + \tfrac{1}{2}\sigma^2(x_*)x_*^2[\bar{\phi}''(x_*) - V''(x_*-)] \\
&= \tfrac{1}{2}\sigma^2(x_*)x_*^2[\bar{\phi}''(x_*) - V''(x_*-)] \\
&= \tfrac{1}{2}\sigma^2(x_*)x_*^2 A\bar{\phi}(x_*).
\end{aligned}$$

Proposition 5.2 follows by combining the last display with (5.3), (5.5) and (5.6).

5.2. *Proof of Proposition* 5.3. We recall the notation $x_\delta = x_* - \delta$, where $\delta = c\sqrt{h} + o(\sqrt{h})$. It follows from the definition (2.1) of $A$ and Taylor's theorem that

$$\begin{aligned}
\triangle \bar{W}_\delta'(x_\delta) &= \bar{\phi}'(x_\delta) - \frac{\bar{\phi}(x_\delta)}{V(x_\delta)}V'(x_\delta) \\
&= [V''(x_*-) - \bar{\phi}''(x_*)]\delta + o(\delta) \\
&= cA\phi(x_*)\sqrt{h} + o(\sqrt{h}).
\end{aligned}$$

As a consequence, the main difficulty in proving (5.2) lies with the term

$$E^x \int_0^{\tau_\delta^h} e^{-rt} \, dL_t^S(x_\delta).$$

As in Section 5.1, we consider the transformation

$$Y_t^{(h)} \doteq \frac{S_{th} - x_\delta}{\sqrt{h}}.$$

Then $Y^{(h)}$ satisfies the SDE

$$dY_t^{(h)} = (\sqrt{h}Y_t^{(h)} + x_\delta)[\sqrt{h}b(\sqrt{h}Y_t^{(h)} + x_\delta) \, dt + \sigma(\sqrt{h}Y_t^{(h)} + x_\delta) \, dB_t].$$

We have the following lemma, whose proof is trivial from the definition of the local time and thus omitted.

LEMMA 5.6. *Suppose $X$ is a semimartingale, and $Y_t \doteq aX_{bt} + v$, where $a > 0, b > 0, v$ are constants. Let $L^Y$ and $L^X$ denote the local times for $Y$ and $X$, respectively. Then, for all $t \geq 0$,*

$$L_t^Y(ax + v) = aL_{bt}^X(x).$$

It follows from the lemma that

$$E^x \int_0^{\tau_\delta^h} e^{-rt} \, dL_t^S(x_\delta) = \sqrt{h} E^x \int_0^{N^h} e^{-rth} \, dL_t^{Y^{(h)}}(0).$$

For any $u \in [0, 1)$, define the process

$$Y_t^* = x_* \sigma(x_*) B_t, \qquad Y_u^* = 0.$$

Also define

$$Q(u) \doteq E L_{u,N}^{Y^*}(0) \qquad \text{where } N \doteq \inf\{n \in \mathbb{N} : Y_n^* \geq 0\}.$$

We have the following result.

PROPOSITION 5.7.

$$\lim_{h \to 0} E^x \int_0^{N^h} e^{-rth} \, dL_t^{Y^{(h)}}(0) = E^x[e^{-r\tau_*}] \int_0^1 Q(u) \, du.$$

Before giving the proof, we show how the desired Proposition 5.3 will follow from Proposition 5.7. We have $E^x e^{-r\tau_*} = V(x)/\bar\phi(x_*)$, and the definitions of $Q$ and $M$, and Lemma 3.1 imply $\int_0^1 Q(u) \, du = x_* \sigma(x_*) \int_0^1 M(u) \, du$. When combined with the expansion given above for $\Delta W_\delta'(x_\delta)$, the left-hand side of (5.2) is equal to

$$hcA\bar\phi(x_*) \frac{V(x)}{\bar\phi(x_*)} x_* \sigma(x_*) \int_0^1 M(u) \, du + o(h),$$

which is exactly the right-hand side of (5.2).

PROOF OF PROPOSITION 5.7. We consider the test function

$$f(x) \doteq \begin{cases} 0, & \text{if } x \leq 0, \\ x, & \text{if } 0 \leq x \leq 1, \\ k, & \text{if } x \geq 2. \end{cases}$$

We require $f(x)$ to be increasing and smooth, except at the point $x = 0$ (the specific choice of $k$ is not important). It follows from the generalized Itô formula and the integration by parts formula that

$$d[e^{-rth} f(Y_t^{(h)})] = -rh e^{-rth} f(Y_t^{(h)}) \, dt + e^{-rth} D^- f(Y_t^{(h)}) \, dY_t^{(h)}$$
$$+ \tfrac{1}{2} e^{-rth} f''(Y_t^{(h)}) \, dY_t^{(h)} \cdot dY_t^{(h)} + e^{-rth} \, dL_t^{X^{(h)}}(0).$$

Without loss of generality, we let $f''(0) = 0$. Now we integrate both sides from 0 to $N^h$ and take expected value.

The first term on the right-hand side will contribute

$$-rhE^x \int_0^{N^h} e^{-rth} f(Y_t^{(h)}) \, dt = -rhE^x \int_{\tau_\delta/h}^{N^h} e^{-rth} f(Y_t^{(h)}) \, dt,$$

since $f(x) = 0$ for $x \leq 0$. We recall the definition (5.4) of $X_{h,u}$. It follows from the strong Markov property that

$$E^x \int_{\tau_\delta/h}^{N^h} e^{-rth} f(Y_t^{(h)}) \, dt \leq k E^x \int_{\tau_\delta/h}^{N^h} \mathbb{1}_{\{Y_t^{(h)} \geq 0\}} \, dt = k E^x G(h, U_h),$$

where

$$U_h \doteq \frac{\tau_\delta}{h} - \left\lfloor \frac{\tau_\delta}{h} \right\rfloor$$

and

$$G(h, u) \doteq E X_{h,u}.$$

By the uniform integrability of $X_{h,u}$ for small $h$ and $u \in [0, 1)$, $E^x G(h, U_h)$ is uniformly bounded for small $h$. Therefore, the expectation of the first term in the right-hand side goes to zero as $h \to 0$.

The second term in the right-hand side contributes (observe that the stochastic integral has expectation zero)

$$\sqrt{h} E^x \int_0^{N^h} e^{-rth} D^- f(Y_t^{(h)}) (\sqrt{h} Y_t^{(h)} + x_\delta) b(\sqrt{h} Y_t^{(h)} + x_\delta) \, dt.$$

Note that the integrand is bounded by $\mathbb{1}_{\{Y_t^{(h)} \geq 0\}}$ up to a proportional constant. It follows exactly as in the case of the first term that the contribution of the second term goes to zero.

The third term in the right-hand side contributes

$$E^x \int_0^{N^h} \tfrac{1}{2} e^{-rth} f''(Y_t^{(h)}) (\sqrt{h} Y_t^{(h)} + x_\delta)^2 \sigma^2 (\sqrt{h} Y_t^{(h)} + x_\delta) \, dt.$$

Since $f''(x) = 0$ for $x < 0$, the expected value equals

$$E^x \int_{\tau_\delta/h}^{N^h} \tfrac{1}{2} e^{-rth} f''(Y_t^{(h)}) (\sqrt{h} Y_t^{(h)} + x_\delta)^2 \sigma^2 (\sqrt{h} Y_t^{(h)} + x_\delta) \, dt.$$

It follows from strong Markov property that the expectation can also be written

$$E^x [e^{-r\tau_\delta} F(h; U_h)],$$

where

$$F(h; u) \doteq E \int_u^{N^h} \tfrac{1}{2} e^{-r(t-u)h} f''(Y_t^{(h)})(\sqrt{h}Y_t^{(h)} + x_\delta)^2 \sigma^2(\sqrt{h}Y_t^{(h)} + x_\delta) \mathbb{1}_{\{Y_t^{(h)} \geq 0\}} \, dt$$

and where $Y^{(h)}$ satisfies the same dynamics with $Y_u^{(h)} = 0$. Since the integrand is bounded due to the fact that $f''(x) = 0$ for all $x \geq 2$, it follows from an analogous argument to the one given in the proof of Lemma 5.5 that:

1. $F(h; u)$ is uniformly bounded for small $h$ and all $u \in [0, 1)$;
2.
$$J(u) \doteq \lim_{h \to 0} F(h; u) = \tfrac{1}{2} E \int_u^N f''(Y_t^*) x_*^2 \sigma^2(x_*) \, dt$$

and the convergence is uniform on any compact subset of $[0, 1)$.

The uniform convergence (on compact sets) of $F$ and Proposition A.3 in the Appendix imply that the expectation of the third term converges to

$$E^x[e^{-r\tau_*}] \int_0^1 J(u) \, du.$$

We omit the details here since an analogous argument is used in the proof of Proposition 5.2.

It remains to calculate the contribution from the term

$$E^x[e^{-r\tau_\delta^h} f(Y_{N^h}^{(h)})] = E^x[e^{-r\tau_\delta} K(h, U_h)],$$

where

$$K(h; u) \doteq E[e^{-r(N^h - u)h} f(Y_{N^h}^{(h)})]$$

with $Y_u^{(h)} = 0$. However, the boundedness and continuity of $f$ ensure the following:

1. $K(h; u)$ is uniformly bounded for all $h$ and all $u \in [0, 1)$.
2.
$$I(u) \doteq \lim_{h \to 0} K(h; u) = E[f(Y_N^*) | Y_u^* = 0]$$

and the convergence is uniform on any compact subset of $[0, 1)$.

Indeed, the first claim is trivial. As for the second claim, let $u^h \to u$. Then as $h \to 0$, $Y^{(h)} \Rightarrow Y^*$. By the Skorohod representation theorem, we can assume $Y^{(h)} \to Y^*$ with probability one, which also implies that $N^h \to N$ with probability one. Therefore, $Y_{N^h}^{(h)} \to Y_N^*$ with probability one. The claim now follows from the dominated convergence theorem. Similarly,

$$E^x[e^{-r\tau_\delta^h} f(Y_{N^h}^{(h)})] \to E^x[e^{-r\tau_*}] \int_0^1 I(u) \, du$$

as $h \to 0$. It is now sufficient to prove

$$I(u) - J(u) = Q(u) \qquad \forall\, u \in [0, 1).$$

This is the same showing

$$E\left[ f(Y_N^*) - \tfrac{1}{2} \int_u^N f''(Y_t^*) x_*^2 \sigma^2(x_*)\, dt - L_{u,N}^{Y^*}(0) \right] = 0,$$

where

$$Y_t^* = x_* \sigma(x_*) W_t, \qquad Y_u^* = 0.$$

But this is a direct consequence from the generalized Itô formula and we complete the proof. □

## APPENDIX

**A.1. Weak convergence of** $(\tau_\delta, U_h)$. For an arbitrary $y > 0$, define the function

$$P^y(x, t) \doteq P\left( \max_{0 \le u \le t} S_u \ge y | S_0 = x \right).$$

We have the following lemma.

LEMMA A.1. *For every fixed $y > 0$, function $P^y \in \mathcal{C}^{1,2}((0, y) \times (0, \infty)) \cap \mathcal{C}((0, y) \times [0, \infty))$ and satisfies the parabolic equation*

$$-\frac{\partial P^y}{\partial t}(x, t) + \mathcal{L} P^y(x, t) = 0, \qquad (x, t) \in (0, y) \times (0, \infty).$$

PROOF. It follows from a standard weak convergence argument that $P^y$ is a continuous function; see, for example, [16]. Let $(x_0, t_0) \in (0, y) \times (0, \infty)$ and define the region

$$D \doteq (x_0 - \varepsilon, x_0 + \varepsilon) \times (t_0 - \varepsilon, t_0).$$

Consider the parabolic equation

$$-\frac{\partial u}{\partial t}(x, t) + \mathcal{L} u(x, t) = 0, \qquad (x, t) \in D,$$

with boundary condition $u = P^y$ on its parabolic boundary. It follows from standard PDE theory that there exists a classical solution $u$ [9]. It remains to show that $u = P^y$ in the domain $D$. Define the stopping time

$$\tau \doteq \inf\{t \ge 0 : (t_0 - t, S_t) \notin D\}.$$

It follows that the process $u(S_t, t_0 - t)$ is a (bounded) martingale. In particular,

$$u(x_0, t_0) = E^{x_0} u(S_\tau, t_0 - \tau) = E^{x_0} P^y(S_\tau, t_0 - \tau) = P^y(x_0, t_0).$$

Here the last equality follows from the strong Markov property. □

For fixed $0 < x < y$, the density of the hitting time $\tau_y$ is defined as

$$p^y(x,t) \doteq \frac{\partial P^y}{\partial t}(x,t).$$

According to the preceding lemma, $p^y$ is continuous in the domain $(0,y) \times (0,\infty)$.

LEMMA A.2.   *Suppose $y_n \to y_*$, then $P^{y_n}(x,t) \to P^{y_*}(x,t)$ and $p^{y_n}(x,t) \to p^{y^*}(x,t)$ uniformly on any compact subset of $(0,y_*) \times (0,\infty)$.*

PROOF.   It suffices to show that $P^{y_n}(x,t) \to P^{y^*}(x,t)$ uniformly on any compact subset. The uniform convergence of $p^{y_n}$ then follows from [9], Section 3.6. Suppose $D \doteq [x_0,x_1] \times [t_0,t_1] \subseteq (0,y_*) \times (0,\infty)$ is a compact subset. In the following, we will denote $P^{y_n}$ and $P^{y_*}$ by $P_n$ and $P$, respectively. Also, we assume $y_n \leq y_*$ for all $n$, which implies that

$$P_n(x,t) \geq P(x,t).$$

An analogous argument can be used for the case $y_n \geq y_*$.

For any $\varepsilon > 0$, we want to show that, for large enough $n$,

$$0 \leq P_n(x,t) - P(x,t) \leq \varepsilon \qquad \forall (x,t) \in D.$$

Define

$$\tau \doteq \inf\{t \geq 0 : S_t \geq y_*\}; \qquad \tau_n \doteq \inf\{t \geq 0 : S_t \geq y_n\}.$$

Since $P$ is continuous, it is uniformly continuous on the compact subset $D$. It follows that there exists a number $h$ such that

$$P(t < \tau \leq t + h | S_0 = x) = P(x,t+h) - P(x,t) \leq \frac{\varepsilon}{2} \qquad \forall (x,t) \in D,$$

and, thus, for all $(x,t) \in D$,

$$P_n(x,t) - P(x,t) = P^x(\tau_n \leq t, \tau > t) \leq P^x(\tau_n \leq t, \tau > t + h) + \frac{\varepsilon}{2}.$$

However, it follows from strong Markov property that, for any $(x,t) \in D$,

$$P(\tau_n \leq t, \tau > t + h | S_0 = x) \leq P\Big(\max_{0 \leq u \leq h} S_u \leq y_* | S_0 = y_n\Big).$$

Note that the right-hand side is independent of $(x,t) \in D$. A proof analogous to that of Lemma 5.5 yields that the right-hand side is dominated by

$$P\Big(\max_{0 \leq t \leq \bar{\sigma}^2 h}\Big[-\frac{c_1}{\underline{\sigma}}t + B_t\Big] \leq \log\frac{y_*}{y_n}\Big).$$

For $n$ big enough, this probability is, at most, $\varepsilon/2$ since $y_n \to y_*$. This completes the proof. □

PROPOSITION A.3.    *Suppose* $f \colon [0, \infty) \to \mathbb{R}$ *is a bounded continuous function with*

$$\lim_{x \to \infty} f(x) = 0,$$

*and* $g \colon [0, 1) \to \mathbb{R}$ *a continuous, bounded function. Then*

$$\lim_{h, \delta \to 0} E^x \left[ f(\tau_\delta) g \left( \frac{\tau_\delta}{h} - \left\lfloor \frac{\tau_\delta}{h} \right\rfloor \right) \right] = E^x f(\tau_*) \cdot \int_0^1 g(u) \, du.$$

*for all* $x \in (0, x_*)$.

PROOF.    Fix $x \in (0, x^*)$. Let $p_\delta$ and $p$ denote the density of $\tau_\delta$ and $\tau_*$, respectively. We can assume that all $x_\delta$ are close to $x_*$, in the sense that $\delta \leq \delta_0$ for some $\delta_0$, and $x < x_\delta$. Since $f(x) \to 0$ as $x \to \infty$, we have

$$E^x \left[ f(\tau_\delta) g \left( \frac{\tau_\delta}{h} - \left\lfloor \frac{\tau_\delta}{h} \right\rfloor \right) \right] = \int_0^\infty f(s) g \left( \frac{s}{h} - \left\lfloor \frac{s}{h} \right\rfloor \right) p_\delta(s) \, ds.$$

For any $\varepsilon > 0$, there exists $0 < a < M < \infty$ such that

$$\int_{[0,a]} f(s) g \left( \frac{s}{h} - \left\lfloor \frac{s}{h} \right\rfloor \right) p_\delta(s) \, ds \leq \|f\|_\infty \cdot \|g\|_\infty P(\tau_\delta \leq a)$$

$$\leq \|f\|_\infty \cdot \|g\|_\infty P(\tau_{\delta_0} \leq a)$$

$$\leq \varepsilon$$

and

$$\int_{[M,\infty]} f(s) g \left( \frac{s}{h} - \left\lfloor \frac{s}{h} \right\rfloor \right) p_\delta(s) \, ds \leq \max_{M \leq x} |f(x)| \cdot \|g\|_\infty \leq \varepsilon.$$

Note that such choices of $(a, M)$ also make the above inequalities hold when $p_\delta$ is replaced by $p$. Also, since $p_\delta \to p$ uniformly on the compact interval $[\varepsilon, M]$, we have

$$\int_a^M f(s) g \left( \frac{s}{h} - \left\lfloor \frac{s}{h} \right\rfloor \right) |p_\delta(s) - p(s)| \, ds \leq \varepsilon$$

for $\delta$ small enough. It remains to show that, for $h$ small enough,

$$\left| \int_a^M f(s) g \left( \frac{s}{h} - \left\lfloor \frac{s}{h} \right\rfloor \right) p(s) \, ds - \int_a^M f(s) p(s) \, ds \cdot \int_0^1 g(u) \, du \right| \leq \varepsilon.$$

We omit the rather straightforward proof, which follows easily from the uniform continuity of $f$, $p$, and $f \cdot p$ on compact intervals.    $\square$

**A.2. Proof of Lemma 3.1.** We first prove the representation for $H(u)$. Consider the continuously differentiable convex function

$$f(x) \doteq \tfrac{1}{2}(x^+)^2 = \begin{cases} 0, & \text{if } x \leq 0, \\ x^2/2, & \text{if } x \geq 0. \end{cases}$$

It follows from the generalized Itô formula [13] that

$$f(W_{N\wedge n}) = \int_u^{N\wedge n} W_t \mathbb{1}_{\{W_t \geq 0\}} \, dW_t + \tfrac{1}{2} \int_u^{N\wedge n} \mathbb{1}_{\{W_t \geq 0\}} \, dt$$

for all integers $n \in \mathbb{N}$. This yields

$$E(W_{N\wedge n}^+)^2 = E \int_u^{N\wedge n} \mathbb{1}_{\{W_t \geq 0\}} \, dt \qquad \forall\, n \in \mathbb{N}.$$

Letting $n \to \infty$, the right-hand side converges to $E \int_u^N \mathbb{1}_{\{W_t \geq 0\}} \, dt$ by the monotone convergence theorem. Since $W_{N\wedge n}^+ \leq W_N$, the result will follow by dominated convergence if one can show that $EW_N^2$ is finite.

To this end, we observe that

$$EW_N^2 = \sum_{n=1}^{\infty} E(W_N^2 \mathbb{1}_{\{N=n\}})$$

$$= \sum_{n=1}^{\infty} \int_0^{\infty} 2x P(W_N \geq x, N = n) \, dx$$

$$= \sum_{n=1}^{\infty} \int_0^{\infty} 2x P(W_N \geq x | N = n) P(N = n) \, dx.$$

However, on the set $\{N = n\}$, the Brownian motion sample path must cross zero during time interval $(n-1, n]$. Let $\Phi$ denote the cumulative distribution function for the standard normal distribution. For every $t \in [0, 1]$, we have the inequality

$$P(W_t \geq x | W_t \geq 0, W_0 = 0) = 2\Phi(-x/\sqrt{t}) \leq 2\Phi(-x) \qquad \forall\, x \geq 0.$$

Then the strong Markov property easily implies that

$$P(W_N \geq x | N = n) = P(W_n \geq x | W_n \geq 0, W_{n-1} < 0, \ldots, W_1 < 0)$$

$$\leq 2\Phi(-x).$$

Since $N$ is finite with probability one, it follows that

$$EW_N^2 \leq \int_0^{\infty} 4x \Phi(-x) \, dx < \infty.$$

It remains to show the representation for $M(u)$. It follows from Tanaka's formula that

$$W_N = W_N^+ = \int_u^N \mathbb{1}_{\{W_t \geq 0\}} \, dW_t + L_{u,N}^W(0).$$

However, since the preceding proof already implies that $E \int_u^N \mathbb{1}_{\{W_t \geq 0\}} \, dt < \infty$ (and, hence, that the stochastic integral has zero mean), we have

$$EW_N = EL_{u,N}^W(0) = M(u).$$

This completes the proof.

**Acknowledgments.** We would like to thank two anonymous referees and an Associate Editor for their careful reading and many helpful suggestions. In particular, the short proof for Remark 2.3 was given by one of the referees.

## REFERENCES

[1] BATHER, J. A. and CHERNOFF, H. (1966). Sequential decisions in the control of a spaceship. *Proc. Fifth Berkeley Symp. Math. Statist. Probab.* **3** 181–207. Univ. California Press, Berkeley. MR224218

[2] BENSOUSSAN, A. (1984). On the theory of option pricing. *Acta Appl. Math.* **2** 139–158. MR748007

[3] BRAUN, M. (1993). *Differential Equations and Their Applications.* Springer, New York. MR1195529

[4] BRENNAN, M. J. and SCHWARTZ, E. S. (1985). Evaluating natural rescource investments. *Journal of Business* **58** 135–157.

[5] BROADIE, M. and GLASSERMAN, P. (1997). Pricing American-style securities using simulation. *J. Econom. Dynam. Control* **21** 1323–1352. MR1470284

[6] DAYANIK, S. and KARATZAS, I. (2003). On the optimal stopping problems for one-dimensional diffusions. *Stochastic Process. Appl.* **107** 173–212. MR1999788

[7] DÉCAMPS, J. P., MARIOTTI, T. and VILLENEUVE, S. (2004). Irreversible investment in alternative projects. Preprint. Available at http://personal.lse.ac.uk/mariotti. MR2082164

[8] DIXIT, A. K. and PINDYCK, R. S. (1994). *Investment Under Uncertainty.* Princeton Univ. Press.

[9] FRIEDMAN, A. (1964). *Partial Differential Equations of Parabolic Type.* Prentice-Hall, Englewood Cliffs, NJ. MR181836

[10] GUO, X. and SHEPP, L. (2001). Some optimal stopping problems with nontrivial boundaries for pricing exotic options. *J. Appl. Probab.* **38** 647–658. MR1860203

[11] JACKA, S. D. (1991). Optimal stopping and the American put. *Math. Finance* **1** 1–14.

[12] KARATZAS, I. (1988). On the pricing of American options. *Appl. Math. Optim.* **17** 37–60. MR908938

[13] KARATZAS, I. and SHREVE, S. E. (1991). *Brownian Motion and Stochastic Calculus.* Springer, New York. MR1121940

[14] KARATZAS, I. and SHREVE, S. E. (1998). *Methods of Mathematical Finance.* Springer, Berlin. MR1640352

[15] KUNITA, H. (1990). *Stochastic Flows and Stochastic Differential Equations.* Cambridge Univ. Press. MR1070361

[16] KUSHNER, H. (1984). *Approximation and Weak Convergence Methods for Random Processes with Applications to Stochastic System Theory.* MIT Press. MR741469

[17] LAMBERTON, D. (1998). Error estimates for the binomial approximation of American put options. *Ann. Appl. Probab.* **8** 206–233. MR1620362

[18] LAMBERTON, D. (2002). Brownian optimal stopping and random walks. *Appl. Math. Optim.* **45** 283–324. MR1885822

[19] MCDONALD, R. and SIEGEL, D. (1986). The value of waiting to invest. *Quarterly J. Econ.* **101** 707–728.

[20] MERTON, R. C. (1990). *Continuous-Time Finance.* Blackwell, Cambridge.

[21] MUSIELA, M. and RUTKOWSKI, M. (1997). *Martingale Methods in Financial Modelling.* Springer, Berlin. MR1474500

[22] MYNENI, R. (1992). The pricing of the American option. *Ann. Appl. Probab.* **2** 1–23. MR1143390

[23] ROGERS, L. C. G. and WILLIAMS, D. (1994). *Diffusions, Markov Processes and Martingales* **1**, **2**. Cambridge Univ. Press. MR1796539

[24] SALMINEN, P. (1985). Optimal stopping for one-dimensional diffusions. *Mathematische Nachrichten* **124** 85–101. MR827892

[25] SHIRYAYEV, A. N. (1978). *Optimal Stopping Rules.* Springer, New York. MR468067

LEFSCHETZ CENTER FOR DYNAMICAL SYSTEMS
BROWN UNIVERSITY
PROVIDENCE, RHODE ISLAND 02912
USA
E-MAIL: dupuis@cfm.brown.edu